\documentclass{amsart}

\usepackage{amsmath}
\usepackage{amsfonts}
\usepackage[utf8]{inputenc}
\usepackage[numbers]{natbib}
\usepackage{graphicx}
\usepackage{mathtools}
\usepackage{amssymb}
\usepackage{amsthm}
\usepackage{tikz-cd} 
\usepackage{enumerate}
\usepackage[english]{babel}
\usepackage{hyperref}
\usepackage{tikz}
\usetikzlibrary{positioning} 
\usepackage{nomencl}

\theoremstyle{definition}
\newtheorem{cor}{Corollary}[section] % same for example numbers
\newtheorem{prop}[cor]{Proposition}
\newtheorem{theorem}[cor]{Theorem}

\newtheorem{conj}[cor]{Conjecture}

\newtheorem{lemma}[cor]{Lemma}

\numberwithin{equation}{section}

\DeclareMathOperator{\su}{SU}

\DeclareMathOperator{\psuu}{PSU}

\DeclareMathOperator{\slg}{SL}

\DeclareMathOperator{\gu}{GU}

\DeclareMathOperator{\syl}{Syl}
\DeclareMathOperator{\gl}{GL}

\DeclareMathOperator{\spdd}{Sp}

\DeclareMathOperator{\franz}{char}

\begin{document}

%%
%% The title of the paper goes here.  Edit to your title.
%%

\title{Proving a conjecture for fusion systems on a class of groups}

%%
%% Now edit the following to give your name and address:
%% 

\author{Patrick Serwene}
\address{Technische Universität Dresden, Faculty of Mathematics, 01062 Dresden, Germany}
\email{patrick.serwene@tu-dresden.de}

%%
%% If there is another author uncomment and edit the following.
%%

%\author{Second Author}
%\address{Department of Mathematics, University of South Carolina,
%Columbia, SC 29208}
%\email{second@math.sc.edu}
%\urladdr{www.math.sc.edu/$\sim$second}

%%
%% If there are three of more authors they are added in the obvious
%% way. 
%%

%%%
%%% The following is for the abstract.  The abstract is optional and
%%% if not used just delete, or comment out, the following.
%%%

\begin{abstract}
We prove the conjecture that exotic and block-exotic fusion systems coincide holds for all fusion systems on exceptional $p$-groups of maximal nilpotency class, where $p \geq 5$. This is done by considering a family of exotic fusion systems discovered by Parker and Stroth. Together with a previous result by the author, which we also generalise in this paper, and a result by Grazian and Parker this implies the conjecture for fusion systems on such groups. Considering small primes, there are no exotic fusion systems on $2$-groups of maximal class and for $p=3$, we prove block-exoticity of two exotic fusion systems described by Diaz--Ruiz--Viruel.
\end{abstract}

%%
%%  LaTeX will not make the title for the paper unless told to do so.
%%  This is done by uncommenting the following.
%%

% \maketitle

%%
%% LaTeX can automatically make a table of contents.  This is done by
%% uncommenting the following:
%%

%\tableofcontents

%%
%%  To enter text is easy.  Just type it.  A blank line starts a new
%%  paragraph. 
%%

\maketitle

\section{Introduction}
Let $p$ be a prime and $P$ a finite $p$-group. A \textit{fusion system} on $P$ is a category whose objects are the subgroups of $P$, we say \textit{category on $P$}, and whose morphisms are injective group homomorphisms between the subgroups of $P$ fulfilling certain assumptions. If these morphisms satisfy two more additional axioms, we call the fusion system \textit{saturated}, but for convenience we drop this term for this paper and thus mean saturated fusion system by fusion system henceforth. All finite groups $G$ determine a fusion system $\mathcal F_P(G)$ on a Sylow $p$-subgroup $P$ of $G$, if we define the morphisms to be the conjugation maps induced by a fixed element in $G$. A fusion system that can be constructed in this way is called \textit{realisable}, whereas a fusion system which is not realisable is called \textit{exotic}. Fix an algebraically closed field $k$ with $\franz k=p$ and let $b$ be a block of $kG$. In this setting, we can define a fusion system on a defect group $P$ of $b$ by again defining the morphisms to be certain conjugation maps induced by an element in $G$. This fusion system is then denoted by $\mathcal F_{(P,e_P)}(G,b)$, where $(P,e_P)$ is a maximal $b$-Brauer pair. Not every fusion system $\mathcal F$ can be constructed in this way, if it can we call $\mathcal F$ \textit{block-realisable}, otherwise \textit{block-exotic}.\\
The following fact is a consequence of Brauer's Third Main Theorem (see \cite[Theorem 3.6]{radha}): If $G$ is a finite group and $b$ is the principal $p$-block of $kG$, i.e. the block corresponding to the trivial character, with maximal $b$-Brauer pair $(P,e_P)$, then $P \in \syl_p(G)$ and $\mathcal{F}_{(P,e_P)}(G,b)=\mathcal{F}_P(G)$. In particular, any realisable fusion system is block-realisable. The converse is still an open problem and has been around in the form of the following conjecture for  a  while, see \cite[Part IV,7.1]{ako} and \cite[9.4]{david}:

\begin{conj}
\label{XX}
\textit{If $\mathcal F$ is an exotic fusion system, then $\mathcal F$ is block-exotic.}
\end{conj}

Some evidence for this conjecture to hold has been collected for example in \cite{solom}, \cite{david}, \cite{ks} or \cite{serwene}. In this paper, we prove this conjecture for all fusion systems on exceptional $p$-groups of maximal class for $p \geq 5$ and for two fusion systems of maximal class on $3$-groups. We recall some of these terms.\\
Let $n \in \mathbb{N}_{\geq 1}$. A $p$-group has \textit{maximal class} if its order is $p^n$ and its nilpotency class $n-1$. Assume $P$ is of maximal class and has order $p^n$, with $n \geq 4$. Define then $\gamma_2(P)=[P,P]$ and $\gamma_{i+1}(P)=[\gamma_i(P),P]$ for $i \geq 2$. A group of the form $C_P(\gamma_i(P)/\gamma_{i+2}(P))$ for $2 \leq i \leq n-2$ is then called a $2$-step centraliser of $P$. If $P$ has more than one $2$-step centraliser, it is called \textit{exceptional}.\\
For $p=2$, there are no exotic fusion systems on groups of maximal class. In fact, it is conjectured that the fusion systems handled in \cite{solom} and \cite{david} are the only exotic fusion systems on $2$-groups. For $p=3$, fusion systems on groups of maximal class have been classified in \cite{drv} and some cases missing in this classification have been added in \cite{parkersem2}. For $p \geq 5$, Grazian--Parker classified all fusion systems on exceptional groups of maximal class in \cite{grazpark} and they found that there are two exotic families: The Parker--Semeraro systems, see \cite{parkersem}, for which block-exoticity has been proven in \cite{serwene} and the Parker--Stroth systems, see \cite{parkerst}, for which we prove block-exoticity here. In particular, this article completes the proof of Conjecture \ref{XX} for all fusion systems on exceptional $p$-groups of maximal class if $p \geq 5$, obtaining the following main theorem:

\begin{theorem}
\label{val}
\textit{Let $p \geq 5$ and $\mathcal F$ be a fusion system on an exceptional $p$-group with maximal class. Then Conjecture \ref{XX} holds for $\mathcal F$}.
\end{theorem}

For $p=3$, if we let $B$ be one of the Blackburn groups $B(3,4;0,2,0)$ or $B(3,2k;0,\gamma,0)$, there has been a reduction of Conjecture \ref{XX} for exotic fusion systems on $B$ to blocks of quasisimple groups in \cite{afaf}. We prove the conjecture for fusion systems on these groups in this article too:

\begin{theorem}
\label{morge}
\textit{Let $\mathcal F$ be an exotic fusion system on one of the groups $B(3,4;0,2,0)$ or $B(3,2k;0,\gamma,0)$ for $k \geq 3$ and $\gamma \neq 0$. Then Conjecture \ref{XX} holds for $\mathcal F$}.
\end{theorem}

See \cite{ako} for details on (block) fusion systems or Section 2 of \cite{serwene} for a more compact overview of the terms needed.\\
In Section \ref{threegroups}, we prove block-exoticity of the exotic fusion systems on the two Blackburn groups and in Section \ref{parker} we do the same for the Parker--Stroth systems. Finally, in Section \ref{rooof}, we combine some results to prove our main theorem.

\section{Fusion Systems on $3$-Groups}
\label{threegroups}
We briefly recall the groups on which the exotic fusion systems we are dealing with in this chapter are defined. See \cite{drv} for details. The non-cyclic $3$-groups of maximal class of order greater than $3^3$ are the groups $B(3,r;\beta,\gamma,\delta)$ of order $3^r$ defined on a set of generators $\{s,s_1,\dots,s_{r-1}\}$. These generators fulfill certain relations, see \cite[Theorem A.2]{drv}, which depend on the parameters $(\beta,\gamma,\delta)$. Theorem 5.9 of \cite{drv} classifies all fusion systems $\mathcal F$ on such groups with $r \geq 4$ and at least one proper $\mathcal F$-essential subgroup. However, it turned out that this classification was incomplete and the fusion systems that have been missed are given in \cite[Theorem 1.1]{parkersem2}.\\
Here we consider exotic fusion systems on the groups $B(3,4;0,2,0)$ and $B(3,2k;0,\gamma,0)$ with $k \geq 3$ and $\gamma \in \{1,2\}$. These exotic fusion systems have been studied with respect to Conjecture \ref{XX} in \cite{afaf} and the following observation has been made:

\begin{theorem}
\label{morge}
\textit{Let $\mathcal F$ be an exotic fusion system on one of the groups $B(3,4;0,2,0)$ or $B(3,2k;0,\gamma,0)$ for $k \geq 3$ and $\gamma \neq 0$. If $\mathcal F$ is block-realisable, it is block-realisable by a block of a quasisimple group with $3'$-centre}.
\end{theorem}

\textit{Proof.} This follows immediately from combining Theorems 1.1, 1.2 and 1.3 from \cite{afaf}. \hfill$\square$\\

This reduction is a key step in proving that these exotic fusion systems are block-exotic as well, which we undertake now.

\begin{theorem}
\label{mennje}
\textit{Let $\mathcal F$ be a fusion system on one of the groups $B(3,4;0,2,0)$ or $B(3,2k;0,\gamma,0)$ for $k \geq 3$ and $\gamma \neq 0$, then Conjecture \ref{XX} holds for $\mathcal F$}.
\end{theorem}

\textit{Proof.} Let $B$ be one of the groups $B(3,4;0,2,0)$ or $B(3,2k;0,\gamma,0)$ for $k \geq 3$ and $\gamma \neq 0$ and $\mathcal F$ an exotic fusion system on $B$. Using the classification in \cite[Theorem 5.9]{drv}, we see that the only proper $\mathcal F$-essential subgroup of $B$ is $C_3 \times C_3$. Further, by Theorem \ref{morge}, if there is a finite group having a block such that its fusion system is one of the aforementioned exotic systems, we can assume that this group is quasisimple with $3'$-centre. Let $G$ be such a group having a block $b$ with maximal $b$-Brauer pair $(B,e)$ and assume $\mathcal F=\mathcal F_{(B,e)}(G,b)$. Further note that by the proof of \cite[Theorem 5.9]{drv}, Part ``Exoticism", the $3$-rank of $G$ must be two. We use the classification of finite simple groups to show that such a group $G$ cannot exist.\\
Firstly, assume $G/Z(G)$ is an alternating group $\mathfrak A_m$. Then $B$ is isomorphic to a Sylow $3$-subgroup of either $\mathfrak S_{3}$ or $\mathfrak S_{6}$. Checking their orders, this is clearly a contradiction.\\
Assume now that $G$ is a finite group of Lie type over a field of characteristic $q$. Firstly, if $q=3$, then $B$ must be a Sylow $3$-subgroup of $G$ by \cite[Theorem 6.18]{caen}, which is not possible since $\mathcal F$ is exotic. Thus, assume $q \neq 3$. In the proof of \cite[Theorem 5.9]{drv}, all finite groups of Lie type over a field of characteristic different from $3$ with $3$-rank two are listed. The possibilites are $L^{\pm}_3(q^l), G_2(q^l),$ $^3D_4(q^l)$ or $^2F_4(q^l)$.\\
Assume first $G=L^{\pm}_3(q^l)$. By \cite[Sections 2.2, 2.3]{hiss}, the defect group of a $3$-block of $kG$ is non-abelian only for the principal $3$-block, which can not give rise to an exotic fusion system. Note that we can exclude the Tits group $^2F_4(2)'$ and its automorphism group $^2F_4(2)$, since $|^2F_4(2)|_3=3^3$. If $G=G_2(q^l)$, the defect groups of $3$-blocks of $kG$ are either Sylow $3$-subgroups or abelian by \cite[2]{hiss}. The same is true for $G$=$^3D_4(q^l)$ by \cite[Proposition 5.4]{derm} and for $G$=$^2F_4(q^l)$ with $q^l \neq 2$ by \cite[Bemerkung 2]{malle}. \\
Finally, assume $G/Z(G)$ is one of the sporadic groups. By, \cite[Theorem 9.22]{david}, the fusion system of a block of a sporadic group can not be exotic. This proves the theorem. \hfill$\square$\\

The $3$-groups handled in this theorem are of maximal nilpotency class. Note that there are more exotic fusion systems on $3$-groups of maximal class by \cite[Theorem B.5]{grazpark}. However, since the reduction in \cite{afaf} is only for these two Blackburn groups, the fusion systems not examined here will require different methods in future work.

\section{The Parker--Stroth Systems}
\label{parker}
Let $p \geq 5$. For context, we briefly recall the Parker--Stroth systems, see \cite{parkerst} for details. Fix $k$ to be a field of characteristic $p$, let $k[x,y]$ be the polynomial algebra in two commuting variables and $V_m$ the $(m+1)$-dimensional subspace of $k[x,y]$ consisting of homogeneous polynomials of degree $m$, where $p-1 \geq  m \geq 1$. Set $L=k^\times \times \gl_2(k) $ and $Q=V_m \times k^+$. Furthermore define $S_0=\{1\} \times \left \{ \begin{psmallmatrix}1 & 0\\ \gamma & 1\end{psmallmatrix} \mid \gamma \in k \right \} \leq L$. Then $L$ acts on $Q$ as follows: for $(v,a) \in Q$ and $(t,A) \in L$ define $(v,a)^{(t,A)}:=(t(v \cdot A), t^2(\det A)^mz)$. Set $P=LQ$ and $S=S_0Q$. The group $K$ consists of all elements in $\slg_3(k)$, where the first row is equal to $(1,0,0)$. Let further $C=\left \{\begin{psmallmatrix} 1&0&0 \\ \alpha & \theta & 0 \\ \delta & \epsilon & \theta^{-1} \end{psmallmatrix} \mid \alpha, \delta, \epsilon \in k, \theta \in k^\times \right \}$. If $p$ and $m$ are odd, it is then shown in \cite[Theorem 2.9]{parkerst} that $P/C_L(Q)$ has a subgroup $W_0$ such that $N_{P/Z(P)}(W_0Z(P)/Z(P)) \cong C$. The group $G$ is then defined to be the amalgamated product $G=P/C_L(Q) \ast_C K$. It should be noted that the group $G$ is infinite and $|S|=p^{p-1}$.
The exotic systems are then defined to be the fusion systems $\mathcal F_S(G)$, see \cite[Lemma 3.4]{parkerst}.\\
For some considerations that follow we also need the following groups $W=O_p(K),$ $P_1=P/C_L(Q)$ and $K_0=N_G(W_0)$. Fix this notation for this section.\\

Before we prove that the Parker--Stroth systems are block-exotic, we need some auxiliary results. Recall that a fusion system on a $p$-group $S$ is called reduction simple if $S$ has no non-trivial proper strongly $\mathcal F$-closed subgroups.

\begin{lemma}
\label{redps}
\textit{Let $\mathcal F$ be a Parker--Stroth system. Then $\mathcal F$ is reduction simple.}
\end{lemma}

In the proof of this statement, we directly use some of the structural information about the groups involved from \cite{parkerst}.\\

\textit{Proof of Lemma \ref{redps}.} Assume $1 \neq N \leq S$ is strongly $\mathcal F$-closed. In particular, $N \unlhd S$ which implies $Z(S) \leq N$. Since we have $W=\langle Z(S)^{K_0} \rangle$, we also get $W \leq N$. Hence $Z_2(S) \leq [Q,S] \leq N$ and finally $S=QW=\langle Z_2(S)^{P_1} \rangle W \leq N$, which implies the claim. \hfill$\square$

\begin{prop}
\cite[Proposition 4.3]{serwene} 
\textit{Let $G$ be a quasisimple finite group and denote the quotient $G/Z(G)$ by $\overline G$. Suppose either $\overline{G}=G(q) \neq E_8(q)$ is a finite group of Lie type and $p$ is a prime number $\geq 5$ or $G(q) = E_8(q)$ and $p>5$. Assume $(p,q)=1$. Let $D$ be a $p$-group such that $Z(D)$ is cyclic of order $p$ and $Z(D) \subseteq [D,D]$. If $D$ is a defect group of a block of $G$, then there are $n, k \in \mathbb{N}$ and a finite group $H'$ with $\slg_n(q^l) \leq H' \leq \gl_n(q^l)$ $($or $\su_n(q^l) \leq H' \leq \gu_n(q^l))$ such that there is a block $c$ of $H'$ with non-abelian defect group $D'$ such that $D'/Z$ is of order $|D/Z(D)|$ for some $Z \leq D' \cap Z(H')$}.
\end{prop}

Note that this proposition was proved for $p \geq 7$ without restrictions on $E_8$. However, since $5$ is a bad prime for $E_8$, we have to consider this group separately in this case.

\begin{theorem}
\label{parkerstroth}
\textit{Let $\mathcal F$ be a Parker--Stroth system defined on $S$. Then $\mathcal F$ is block-exotic}.
\end{theorem}

\textit{Proof.} Let $H$ be a finite group with a block $b$ having maximal $b$-Brauer pair $(S,e_S)$ such that $\mathcal F=\mathcal F_{(S,e_S)}(H,b)$. Since $Z(S)$ is cyclic, we can apply \cite[Theorem 3.5]{serwene} to see that we can assume that $H$ is quasisimple. We use the classification of finite simple groups to exclude all possibilites for $H$.\\
Firstly, assume $H$ is an alternating group $\mathfrak A_m$. Then by \cite[Theorem 6.2.45]{jk} $S$ is isomorphic to a Sylow $p$-subgroup of some symmetric group $\mathfrak S_{pw}$ with $w \leq p-1$. Define the cycle $\sigma_i=((i-1)p+1,\dots,ip)$ and the subgroup $S'=\langle \sigma_1, \dots, \sigma_{p-1} \rangle\leq \mathfrak A_m$. Then $S' \in \syl_p(\mathfrak A_m)$. But this group is abelian, so $S$ cannot be the defect group of a block of $\mathfrak A_m$.\\
Next, assume $H$ is a group of Lie type. First assume the latter group is defined over a field of characteristic $p$, then $S$ is a Sylow $p$-subgroup of $H$ by \cite[Theorem 6.18]{caen} and thus $\mathcal F$ cannot be exotic. In particular, we can assume $H$ is defined over a field of order coprime to $p$. Recall that we have $|S|=p^{p-1}$ and $|Z(S)|=p$.\\
First assume either $G \neq E_8(q)$ or $G=E_8(q)$ and $p>5$. Let $H'$, $D'$ be as in the assertion of the previous proposition. Assume first $H' \leq \gl_n(q^l)$ and let $a$ be such that $|q^l-1|_p=p^a$. Then, since $\slg_n(q^l) \leq H'$, we have $|D'|=|S/Z(S)| \cdot |Z|=p^{p-2}|Z| \leq p^{p-2}|Z(H')| \leq p^{p-2}|Z(\slg_n(q^l))|\leq p^{p+a-2}$. Now the block of $k\gl_n(q^l)$ covering $c$ has a defect group of order at most $p^{p+2a-2}$. But it is a well-known fact, that (non-abelian) defect groups of $\gl_n(q^l)$ have order at least $p^{pa+1}$, see \cite[Theorem 3C]{torquay}. Thus, $p^{pa+1} \leq p^{p+2a-2}$, which is a contradiction if $a \geq 1$, which we can assume in a non-abelian case. The case $H' \leq \gu_n(q^l)$ can be shown in the same fashion by considering the $p$-part of $q^l+1$ instead of $q^l-1$.\\
Now assume $G=E_8(q)$ and $p=5$. By \cite[Theorem 7.7]{bdr}, we can reduce to quasi-isolated blocks. The quasi-isolated $5$-blocks of $G$ are studied in \cite[Section 6.3]{kessarmalle}. There, in Tables 7 and 8 and in the respective supplement in Tables 1 and 2, structural information about the possible defect groups are given. By \cite[Proposition 6.10]{kessarmalle}, the defect groups of these blocks are abelian if the respective relative Weyl group is a $5'$-group. Checking first Table 8, we see that all relative Weyl groups there have $5'$-order and thus all blocks handled in Table 8 have an abelian defect group, so we can exclude all these blocks. Looking at Table 7, also by \cite[Proposition 6.10]{kessarmalle}, we can first restrict ourselves to the blocks described in lines 1, 3, 7, 9, 10, 11, 13, 14, 16 and 19. Looking at the Weyl groups, we can immediately exclude lines 9 and 14, because those blocks have abelian defect groups. Checking the orders of the remaining defect groups, we see that all their orders are bigger than $5^4$, which means that none of these blocks can have $S$ as defect group. For the cases covered in the supplement, all the Weyl groups in lines 46--51 are $5'$ and thus the corresponding blocks are of abelian defect. This only leaves us with line 45, but here we see that the defect group has order $5^9$ and thus can also not be equal to $S$.\\
Finally, assume $H$ is one of the sporadic groups. By \cite[Theorem 9.22]{david}, the fusion system of a block of a sporadic group can not be exotic. This proves the theorem.
\hfill$\square$

\section{Proof of Main Theorem}
\label{rooof}
In this section, we combine the main theorem of the previous section with previous restults to prove Theorem \ref{val}. For this, we generalise a previous result by the author, which has already appeared in this more general form in \cite{phd}.

\begin{theorem}	
\label{psx}
\textit{Conjecture $\ref{XX}$ is true for all fusion systems $\mathcal F$ on a Sylow $p$-subgroup of $G_2(p^n)$ or $\psuu_4(p^n)$ for all primes $p$ and $n \in \mathbb N$}.
\end{theorem}

\textit{Proof.} By \cite[Main Theorem]{vanbeek}, the only exotic fusion systems on such groups are the 27 exotic Parker--Semeraro systems. For these, the conjecture has been proven in \cite[Theorem 1.2]{serwene}. \hfill$\square$\\

\textit{Proof of Theorem \ref{val}}. By \cite[Theorem B]{grazpark}, such an exotic fusion system is either one of the 27 Parker--Semeraro systems or a Parker--Stroth system. The Parker--Semeraro systems are block-exotic by Theorem \ref{psx}, and the Parker--Stroth systems by Theorem \ref{parkerstroth}. Hence the theorem follows. \hfill$\square$

\section*{Acknowledgements}
The author would like to thank the Isaac Newton Institute for Mathematical Sciences, Cambridge, for support and hospitality during the programme Groups, Representations and Applications: new perspectives (GRA2) where work on this paper was undertaken. This work was supported by EPSRC grant no EP/R014604/1. In particular, the author would like to thank Valentina Grazian for very interesting discussions and suggestions during the aforementioned program.\\
Also, the author would like to thank Radha Kessar for further suggestions in pursuing this project and the University of Cumbria, Ambleside, for hosting the workshop ``Structure of Group Algebras over Local Rings", funded by the EPSRC through grants EP/T004592/1 and EP/T004606/1, Representation Theory over Local Rings, which provided a platform for these discussions.

\end{document}